\documentclass[11pt]{amsart}
\usepackage{amsmath,amssymb,epsfig,color,latexsym,amsthm,subfig, mathrsfs}
\usepackage{wrapfig,lipsum,floatrow,sidecap,diagbox,amsfonts,comment}
\usepackage{graphicx}
\usepackage{multirow}
\usepackage{hhline}
\usepackage{url}
\usepackage{epstopdf}
\usepackage{mathtools}
\usepackage{pgfplots}
\usepackage{booktabs}
\graphicspath{ {images/} }
\pgfplotsset{width=8cm,compat=1.9}
\allowdisplaybreaks
\def\p{\prime}
\def\L{\mathcal{L}}

\usepackage{longtable}
\newtheorem{theorem}{Theorem}[section]
\newtheorem{conjecture}{Conjecture}[section]
\newtheorem{corollary}{Corollary}[section]

\newtheorem{remark}{Remark}[section]

\begin{document}

\title[The average genus of oriented rational links]{The average genus of oriented rational links with a given crossing number}
\author{ Dawn Ray$^{\ddag}$, Yuanan Diao$^{\ast}$$^{\ddag}$}
\address{
$^{\ddag}$ Department of Mathematics and Statistics\\
University of North Carolina at Charlotte\\
 Charlotte, NC 28223, USA\\}
 \email{dmray@uncc.edu}

\begin{abstract}
In this paper, we enumerate the number of oriented rational knots and the number of oriented rational links with any given crossing number and minimum genus. This allows us to obtain a precise formula for the average minimal genus of oriented rational knots and links with any given crossing number. 
\end{abstract}

\maketitle
\section{Introduction}
A Seifert surface for an oriented rational link type $\L$ is an orientable surface whose boundary is a representation of $\L$. The minimal genus of $\L$, denoted by $g(\L)$, is defined as the minimum over the genera of all Seifert surfaces of $\L$. In the case that $\L$ is an alternating link with a reduced alternating diagram $D$, it is known that $g(D)=g(\L)$, where $g(D)$ is the genus of the Seifert surface obtained from $D$ using the Seifert's algorithm \cite{Gabai1986}. Furthermore, $g(D)$ can be computed using the formula $2g(D)=c(D)-s(D)-\mu(D)+2$, where $c(D)=Cr(\L)$ is the crossing number of $D$ (which equals the minimum crossing number $Cr(\L)$ of $\L$),  $s(D)$ is the number of Seifert circles in $D$ and $\mu(D)=\mu(\L)$ is the number of components in $D$ ($\L$) \cite{A,BZ}. This formula allows us to compute $g(\L)$, a knot invariant of $\L$, from any reduced alternating diagram of $\L$. In this paper, we shall focus on the rational knots and links. 

\begin{figure}[!hbt]
\begin{center}
\includegraphics[scale=.35]{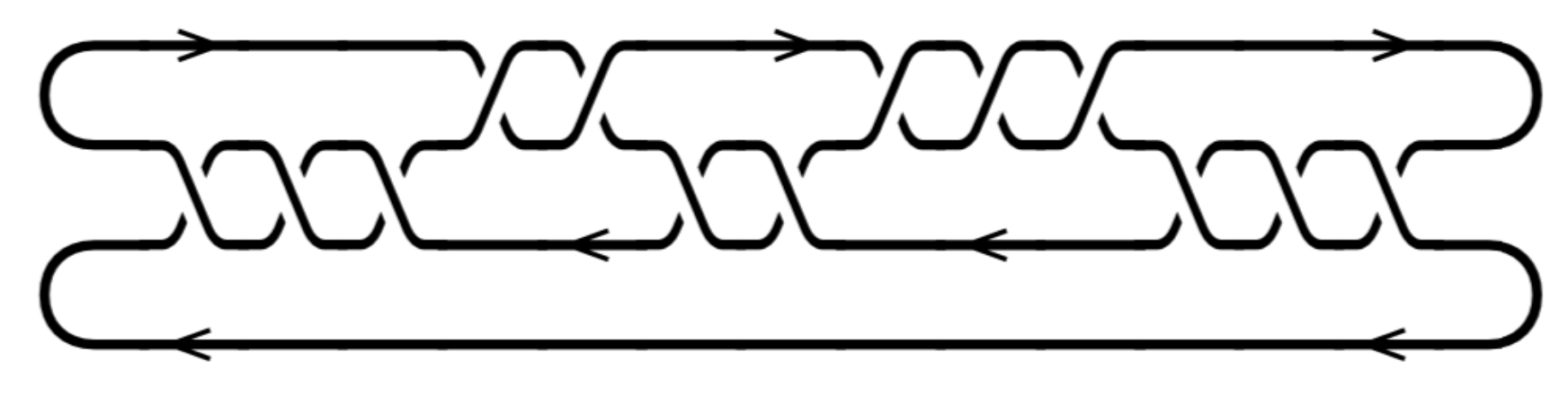}
\end{center}
\caption{A 4-plat based on the vector $[3,2,2,3,3]$.
\label{Figure1}}
\end{figure}

\medskip 
Let $p$ and $q$ be two positive integers with $\rm{gcd}(p,q)=1$ and $0<p<q$. It is known that there exists a unique vector of odd length $(a_1,a_2,...,a_{2m+1})$ with positive integer entries such that 
$$
\displaystyle\frac{p}{q}=\displaystyle\frac{1}{a_1+\frac{1}{a_2+\frac{1}{....\frac{1}{a_{2m}+\frac{1}{a_{2m+1}}}}}}.
$$
We shall present the above in vector form as $p/q=[a_1,a_2,...,a_{2m+1}]$. Given such an odd length vector $[a_1,a_2,...,a_{2m+1}]$ with positive integer entries, one can draw a link diagram called a 4-plat based on this vector as shown in Figure \ref{Figure1}, using the example $56/191=[3,2,2,3,3]$. A rational link is one that has a diagram in the form of a 4-plat. That is, for each rational link, there exist two positive integers $p$ and $q$ such that $\rm{gcd}(p,q)=1$, $0<p<q$, and a 4-plat based on $p/q=[a_1,a_2,...,a_{2m+1}]$ is a link diagram of the link. Furthermore, the minimum crossing number of the rational link is given by $\sum_{1\le j\le 2m+1}a_j$.
A rational link may have up to two components. Throughout the rest of this paper, we shall use the term {\em rational knot} for a rational link with one component, and use the term {\em rational link} when it has two components.

\medskip
The rational knots/links are often used as candidates for testing a random knot generator (to detect for biases for example). The average genus of  rational knots with a given crossing number has recently been studied in \cite{Cohen}. In this paper, we will derive a precise formula for the average genus, not only for the rational knots, but also for the oriented rational links, of any given crossing number.

\section{Oriented rational links and their Seifert circle decompositions}

As we mentioned in the introduction, every oriented rational link can be presented as a 4-plat, which is a reduced alternating diagram, as shown in Figure \ref{Figure1}. We can calculate its genus by finding the number of Seifert circles in its Seifert circle decomposition. It would be more convenient if we can work with the Seifert circle decompositions directly. In order to do so, we need to establish some correspondence relationship between the set of oriented rational links with a given crossing numbers and a set of suitably chosen Seifert circle decompositions. The tools and results presented in this section are based on the work of \cite{Diao2022}.

\medskip
It has been shown in \cite{Diao2022} that every oriented rational link can be represented by a 4-plat corresponding to some rational number $p/q=[a_1,a_2,...,a_{2k+1}]$ such that the long arc at the bottom is oriented from right to left, the leftmost string of horizontal crossings of the 4-plat correspond to $a_1$ and the strand connected to the long arc at the bottom being the under strand. Such a 4-plat is said to be in a {\em preferred standard form}, or just a PS form for short. Let $[b_1, b_2, ..., b_{2k+1}]$ be the signed vector obtained from $[a_1,a_2,...,a_{2k+1}]$ such that $b_i = \pm a_i$ with the sign being the crossing sign of the crossings corresponding to $a_i$ in the 4-plat. For example the signed vector of the 4-plat in Figure \ref{Figure1} is $[-3,-2,-2,-3,-3]$ since all the crossings there have negative signs. 

\begin{figure}[htb!]
\begin{center}
\includegraphics[width=4.4in,height=0.5in]{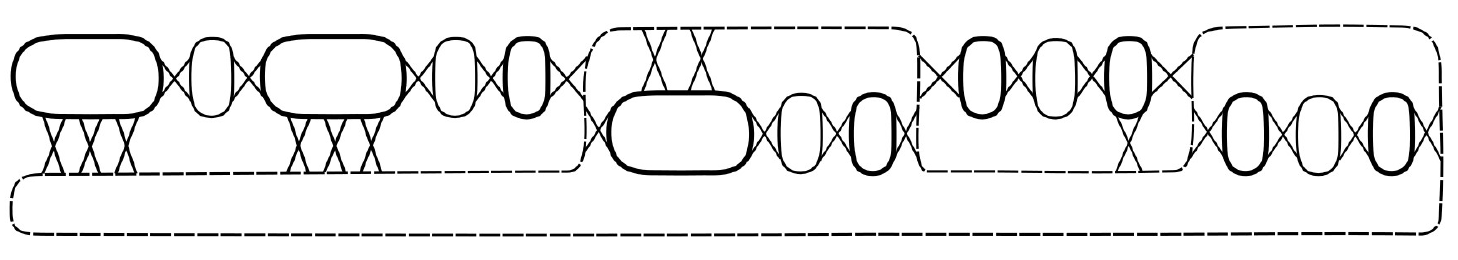}\\
\vspace{0.2in}
\includegraphics[width=5.4in,height=0.5in]{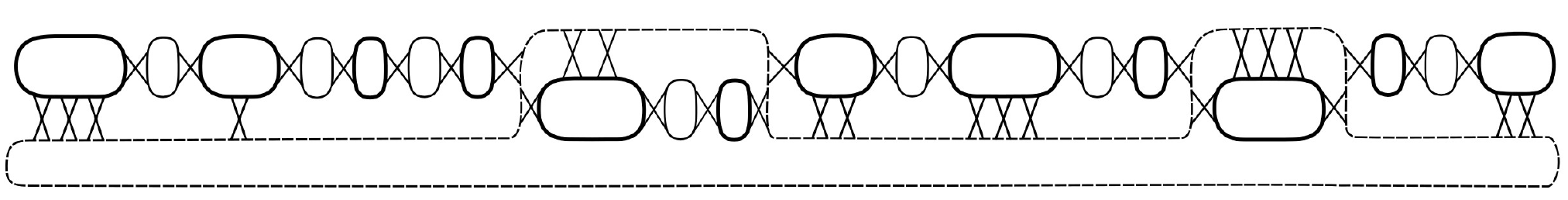}
\end{center}
\caption{Top: A type I R-composition with 5 small and 9 medium Seifert circles; Bottom: A type III R-composition with 7 small and 12 medium Seifert circles. The medium Seifert circles are drawn with thicker lines. Notice that there is a medium Seifert circle in the bottom figure that does not share any crossing with $C$ in this figure, but it can share a crossing with $C$.}
\label{Figure2}
\end{figure}

\medskip
The Seifert circle decomposition of a 4-plat has a special Seifert circle that runs through the entire diagram, namely the one that contains the long arc at the bottom of the 4-plat. Figure \ref{Figure2} shows two typical Seifert circle decompositions of rational knots and links where this special Seifert circle runs through the entire graph. 
For the sake of convenience, we call this Seifert circle the {\em large} Seifert circle and denote it by $C$.
The other Seifert circles may or may not share crossings with $C$. We shall call Seifert circles that CAN share a crossing with $C$ {\em medium} Seifert circles and Seifert circles that cannot share crossings with $C$ {\em small} Seifert circles. We caution that this definition is slightly different from the one in \cite{Diao2022}. For 4-plats that are in PS form, the reader can easily verify that if a small or medium Seifert circle is outside (inside) of $C$, then it can only share positive (negative) crossings with other Seifert circles (including $C$). Furthermore, medium and small Seifert circles that share consecutive crossings outside (inside) of $C$ appear in an alternating order due to the orientations of the circles. 
In \cite{Diao2022}, 4-plats (in PS forms) and their corresponding Seifert circle decompositions are divided into the following 4 types based on the signed vectors $[b_1, b_2, ..., b_{2k+1}]$ of the 4-plats:

\medskip
\begin{center}
(I) $b_1>0$, $b_{2k+1}<0$; (II) $b_1<0$, $b_{2k+1}>0$; (III) $b_1>0$, $b_{2k+1}>0$; (IV) $b_1<0$, $b_{2k+1}<0$.
\end{center}

\medskip
Third, if we rotate a 4-plat in PS form by 180 degrees around an axis that is perpendicular to the horizontal line of the long arc of the 4-plat then change the orientations of the components in the resulting diagram, then we obtain a new 4-plat that is still in the PS form, which is called the {\em reversal} of the original 4-plat. The reversal of a 4-plat is equivalent to its inverse, hence is equivalent to itself since it is known that rational links are equivalent to their inverses. We have the following theorem.

\medskip
\begin{theorem}\label{T1}{\em \cite[Theorem 2.3]{Diao2022}}
Every oriented rational link $\mathcal{L}$ can be uniquely, up to a reversal, represented by a 4-plat in the PS form. 
\end{theorem}

\medskip
Following \cite{Diao2022}, we shall call a Seifert circle decomposition obtained from a 4-plat in the PS form an {\em R-decomposition}. Furthermore, in an R-decomposition we require that the crossing information be retained. Therefore, from an R-decomposition we can recover the corresponding 4-plat in the PS form. Denote by $R_n^I$ and $R_n^{III}$ the sets of type I and III R-decompositions of 4-plats with $n$ crossings respectively. Furthermore, denote by $RS_n^{III}$ the set of type III R-decompositions that are symmetric with respect to the reversal operation. Figure \ref{Figure2} shows two typical type I and III R-decompositions.

\medskip
Let $R_n$ be the set $R_n^I\cup R_n^{III}\uplus RS_n^{III}$, where $\uplus$  means $R_n$ contains an additional copy of $RS_n^{III}$. Let $\Lambda_n$ be the set of all oriented rational knots and links with crossing number $n$.  
We have the following theorem, which is a consequence of \cite[Remark 3.4]{Diao2022}. 

\medskip
\begin{theorem}\label{T2}
There exists a one-to-one and onto mapping $\phi: \Lambda_n\longrightarrow R_n$ such that for each $\L\in \Lambda_n$, there exists a 4-plat $D$ of $\L$ (not necessarily in the PS form) such that  the Seifert circle decomposition of $D$ is the same as $\phi(\L)$.
\end{theorem}

\medskip
Theorem \ref{T2} allows us to work with the set $R_n$ instead of $\Lambda_n$ to extract the information we need in order to compute the genus of any rational knot/link, since the structure of the R-decompositions makes an accurate counting possible as we shall demonstrate in the next section. For the sake of convenience, we shall use an R-decomposition and its corresponding 4-plat via the mapping in Theorem \ref{T2} interchangeably from this point on. 

\medskip
Let $R_n(s)$ be the subset of $R_n$ such that each R-decomposition in $R_n(s)$ contains $s$ Seifert circles. Similarly define $R_n^{I}(s)$, $R_n^{III}(s)$ and $RS_n^{III}(s)$. Then $R_n$ is a disjoint union of the $R_n(s)$'s and $R_n(s)= R_n^I(s)\cup R_n^{III}(s)\uplus RS_n^{III}(s)$. 

\section{The re-construction of R-decompositions from R-templates}

\begin{figure}[htb!]
\begin{center}
\includegraphics[scale=1.0]{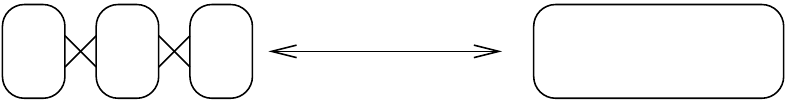}\\
\end{center}
\caption{The reduction and the insertion operation.}
\label{Figure3}
\end{figure}

Consider an R-decomposition with at least one small Seifert circle. If we delete one small Seifert circle and the two crossings it shares with its neighboring Seifert circles, then merge its two neighboring Seifert circles, we obtain a new R-decomposition which has 2 less crossings and 2 less Seifert circles than the original R-decomposition. We call this operation a {\em reduction operation}. The reversed operation of a reduction operation is called an {\em insertion operation}, which is to split a Seifert circle (either small or medium) into two, then insert a small Seifert circle between these two Seifert circles such that the newly inserted small Seifert circle shares a crossing each with its neighboring Seifert circles. This is illustrated in Figure \ref{Figure3}.

\begin{figure}[!ht]
\begin{center}
 \includegraphics[width=3.4in,height=0.7in]{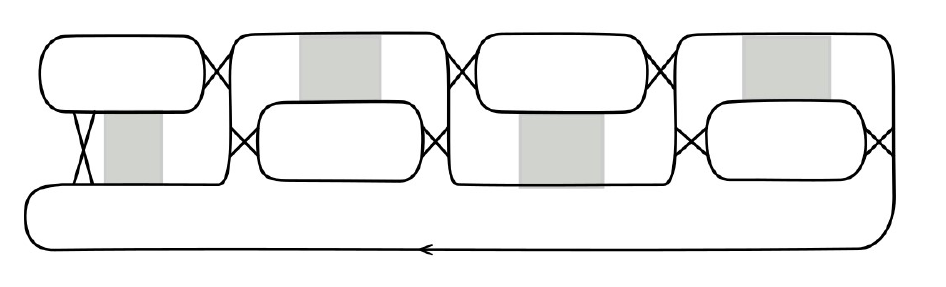}
\end{center}
\caption{The R-template obtained from the R-decomposition in the top of Figure \ref{Figure2} after performing 5 reduction and 7 deletion operations. The crossings deleted by the deletion operations are in the shaded boxes. The number of deleted crossings in each box (from left to right) is: 5, 2, 0, 0. The gray boxes indicate where crossings can be added via the addition operation between the large Seifert circle and the medium Seifert circles.
\label{Figure4}}
\end{figure}

\medskip
Now consider an R-decomposition that contains no small Seifert circles. The operation that deletes a crossing between a medium Seifert circle and the large Seifert circle corresponds to smoothing the crossing in the corresponding link diagram. We call this operation a {\em deletion operation}. The reserved operation of a deletion operation is called an {\em addition operation}. If we delete all but the leftmost and the rightmost crossing for each medium Seifert circle that shares more than two crossings with the large Seifert circle, we arrive at a special R-decomposition satisfying the following two conditions: (i) it contains no small Seifert circles and (ii) each medium Seifert circle in it shares exactly two crossings with the large Seifert circle. We shall call such an R-decomposition an {\em R-template}. By performing reduction operations repeatedly until there are no small Seifert circles left, then perform deletion operations as needed, we can change every R-decomposition to an R-template. Furthermore, it can be shown that this resulting R-template is uniquely determined by the initial R-decomposition and is independent of the order of the reduction operations and deletion operations.
Reserving the operations, we see that every R-decomposition can be reconstructed from an R-template by performing a sequence of addition operations, followed by a sequence of insertion operations. For example, 
there are 15 Seifert circles in the R-decomposition shown at the top of Figure \ref{Figure2}. Of these 5 are small and 9 are medium.  The corresponding diagram has 25 crossings. We can perform 5 reduction operations, which leads to an R-decomposition with 5 Seifert circles and 15 crossings. We can then perform deletion operations to obtain an R-template with 4 medium Seifert circles and 8 crossings as shown in Figure \ref{Figure4}.

\begin{remark}\label{Remark1}{\em
The addition operation and the insertion operation do not change the type of an R-decomposition. Since a type I R-template has an odd number of Seifert circles and a type III R-template has an even number of Seifert circles, when we restrict ourselves to only type I and III R-decompositions (which is what we will do for the rest of this paper), the number of Seifert circles in the R-decomposition determines its type.}
\end{remark}

\begin{figure}[!ht]
\begin{center}
 \includegraphics[width=3.4in,height=0.7in]{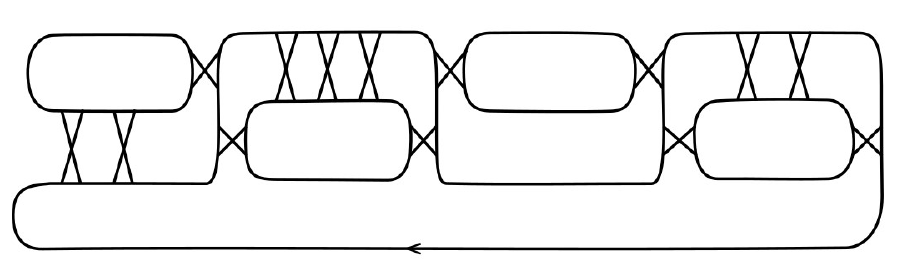}
\end{center}
\caption{An R-decomposition obtained from the R-template in Figure \ref{Figure4} with 6 free crossings added
as shown.
\label{Figure5}}
\end{figure}

\medskip
Let us  denote a type I or III R-decomposition with $n$ crossings, $s$ Seifert circles whose R-template has $q\ge 2$ Seifert circles by $r_n(s,q)$. The following is our adopted procedure for reconstructing $r_n(s,q)$ from a type I or III R-template with $q$ Seifert circles. The template already has a total of $q$ Seifert circles and $2q-2$ crossings, so it is necessary that $s-q\ge 0$ is even (since the rest of the Seifert circles can only be added by pairs via the insertion operation), and we will have to perform $(s-q)/2$ insertion operations  so the resulting R-decomposition will have $s$ Seifert circles as desired. This will increase the number of crossings from $2q-2$ to $2q-2+s-q=s+q-2$. Hence we must have $s+q-2\le n$ and will have to perform $n-s-q+2$ addition operations. We shall call the $2q-2$ crossings in the R-template {\em essential crossings}, and the $n-s-q+2$ crossings to be added via the addition operation  {\em free crossings}. Figure \ref{Figure5} shows an R-decomposition obtained from the R-template in Figure \ref{Figure4} with 6 free crossings added.

\medskip

\section{The enumeration of $|R_n^I(s)\cup R_n^{III}(s)|$}

Recall that $R_n^I(s)\cup R_n^{III}(s)$ is the set of all type I or type III R-decompositions with $n$ crossings and $s$ Seifert circles.  By Remark \ref{Remark1}, if $s$ is odd (even), then the corresponding R-decomposition is of type I (type III).  In this section we prove the following precise formula for $|R_n^I(s)\cup R_n^{III}(s)|$. 

\medskip
\begin{theorem}\label{T4}
\begin{equation}\label{eq_Rn}
|R_n^I(s)\cup R_n^{III}(s)|=
\sum_{j=1}^{\min\{\lceil\frac{s-1}{2}\rceil, \lfloor\frac{n}{2}\rfloor-\lfloor\frac{s-1}{2}\rfloor\}}{{j+\lfloor\frac{s-1}{2}\rfloor-1}\choose{2j-\frac{(-1)^s+1}{2}-1}}{ {n-j-\lfloor\frac{s+1}{2}\rfloor}\choose{j+\lfloor\frac{s-1}{2}\rfloor-1}}.
\end{equation}
\end{theorem}

\medskip
\begin{proof}
Consider the case when $s=2k+1$ for some $k\ge 1$. Consider an R-template $T_{2j+1}$ with $2j+1$, $1\le j\le k$, Seifert circles. $T_{2j+1}$ contains $4j$ crossings. Then the remaining $s-(2j+1)=2k-2j$ Seifert circles are to be inserted via  $k-j$ insertion operations and each such operation adds a small Seifert circle and a medium Seifert circle, as well as two crossings. These insertions can be placed between the large Seifert circle and any of the $2j$ medium Seifert circles in the template as shown in \ref{Figure7} for the simple case $s=5$ and $2j+1=3$. There are ${j+k-1}\choose{2j-1}$ ways to perform the insertions, which lead to a set $TR_{2j+1}$  that contains ${j+k-1}\choose{2j-1}$ distinct R-decompositions obtained from $T_{2j+1}$. Each R-decomposition in $TR_{2j+1}$ contains $j+k$ medium Seifert circles and $2j+2k\le n$ crossings (none of which is a free crossing). The $n-2j-2k\ge 0$ free crossings can be placed between any medium Seifert circle and $C$. We need to point out that $j\le k$ and $j\le \frac{n}{2}-k$ hence we have $1\le j\le \min\{k, \lfloor\frac{n}{2}\rfloor-k\}$.

\begin{figure}[!hbt]
\begin{center}
\includegraphics[width=2.1in,height=0.7in]{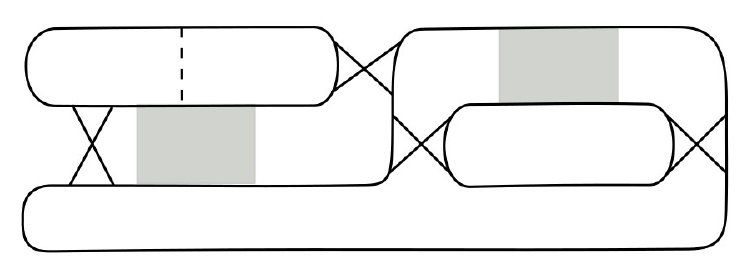}\qquad\includegraphics[width=2.1in,height=0.7in]{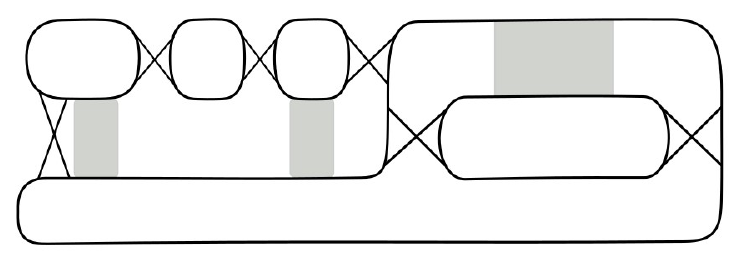}
\end{center}
\caption{Left: An R-template with 3 Seifert circles. The dashed line indicates where an insertion is to be made. Right: The insertion operation splits the medium Seifert circle to two, and adds a small one between them.
\label{Figure7}}
\end{figure}

It follows that the total number of R-decompositions in $R_n^I(s)\cup R_n^{III}(s)$ that can be reconstructed from $T_{2j+1}$ is
${{j+k-1}\choose{2j-1}}{ {n-j-k-1}\choose{j+k-1}}$. Summing over $j$ under the condition $1\le j\le k$, we obtain a precise counting formula for $|R_n^I(s)\cup R_n^{III}(s)|$ when $n$ and $s$ are both odd:
$$
|R_n^I(s)\cup R_n^{III}(s)|=\sum_{j=1}^{\min\{k, \lfloor\frac{n}{2}\rfloor-k\}}{{j+k-1}\choose{2j-1}}{ {n-j-k-1}\choose{j+k-1}}.
$$
Since $k=\frac{s-1}{2}=\lfloor\frac{s-1}{2}\rfloor=\lceil\frac{s-1}{2}\rceil$, $k+1=\frac{s+1}{2}=\lfloor\frac{s+1}{2}\rfloor$ and $\frac{(-1)^s+1}{2}=0$, (\ref{eq_Rn}) holds.

\medskip
The case when $s$ is even can be similarly proven and is left to the reader.
\end{proof}

For example, $|R_7(5)|={{2}\choose{1}}{{3}\choose{2}}=6$,  $|R_8(5)|={{2}\choose{1}}{{4}\choose{2}}+{{3}\choose{3}}{{3}\choose{3}}=13$, $|R_9(5)|={{2}\choose{1}}{{5}\choose{2}}+{{3}\choose{3}}{{4}\choose{3}}=24$, $|R_7(4)|={{1}\choose{0}}{{4}\choose{1}}+{{2}\choose{2}}{{3}\choose{1}}=7$, while $|R_8(4)|={{1}\choose{0}}{{5}\choose{1}}+{{2}\choose{2}}{{4}\choose{1}}=11$.

\section{The enumeration of $|RS^{III}_n(s)|$}

Recall that  $RS^{III}_n(s)$ is the set of all symmetric R-decompositions with $n$ crossings and $s$ Seifert circles. It is necessary that such R-decompositions are of type III.
In this section we prove the following precise formula for $|RS^{III}_n(s)|$. 

\medskip
\begin{theorem}\label{T5}
$|RS^{III}_n(s)|$ is given by the following formula
\begin{equation}\label{RSns}
\left(\frac{(-1)^s+1}{2}\right)\sum_{j=1}^{\min\{\lfloor\frac{s}{2}\rfloor,\lfloor\frac{n}{2}\rfloor+1-\lfloor\frac{s}{2}\rfloor\}}\left(\frac{(-1)^{(j+\lfloor\frac{s}{2}\rfloor)n}+1}{2}\right){{\lfloor{\frac{j+\lfloor\frac{s}{2}\rfloor}{2}}\rfloor-1}\choose{j-1}}{{\lfloor\frac{n}{2}\rfloor-\lceil \frac{j+\lfloor\frac{s}{2}\rfloor}{2}\rceil}\choose{\lfloor{\frac{j+\lfloor\frac{s}{2}\rfloor}{2}}\rfloor-1}}.\\
\end{equation}
\end{theorem}

\medskip
\begin{proof}
Recall that an odd value of $s$ indicates that the R-decomposition is of type I hence cannot be symmetric. Thus if $s=2k+1$ for some positive integer $k$, then $|RS^{III}_n(s)|=0$. The fraction $\frac{(-1)^s+1}{2}$ ensures this. Now we will concentrate on the case  $s=2k$ for some positive integer $k$ (so $k=\frac{s}{2}=\lfloor\frac{s}{2}\rfloor$), we will consider two separate cases: $n$ is even and $n$ is odd.

\medskip
Consider first the case when $n$ is even. An R-template must also have an even number $2j$ of Seifert circles for some $1\le j\le k$. The remaining $2k-2j$ Seifert circles are to be inserted via  $k-j$ insertions  in a symmetric manner. 
If $k-j$ is odd, then one insertion must be placed in the middle of the template which adds one more medium Seifert circle to the $2j-1$ medium Seifert circles in the template. The remaining $k-j-1$ insertions must be distributed symmetrically among the $2j$ medium Seifert circles. Figure \ref{Figure8} shows the case of $s=2k=10$ and $2j=4$.
\begin{figure}[!hbt]
\begin{center}
\includegraphics[scale=0.6]{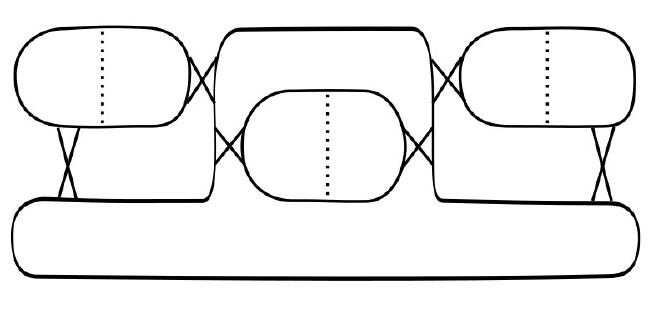} \qquad \includegraphics[scale=0.6]{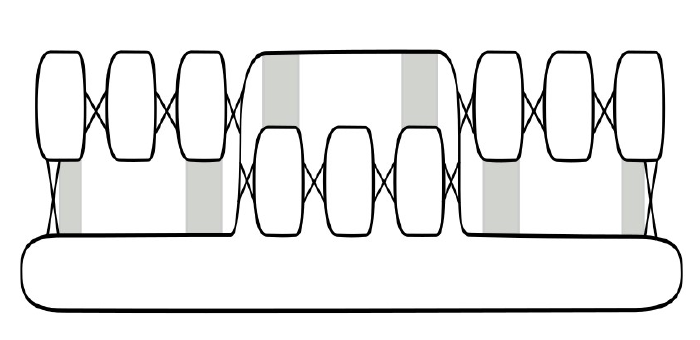}
\end{center}
\caption{A symmetric insertion distribution for an R-decomposition with 10 Seifert circle (right) built on an R-template with 4 Seifert circles (left).
\label{Figure8}}
\end{figure}

The total number of such insertions is ${{\frac{k-j-1}{2}+j-1}\choose{j-1}}={{\frac{k+j-1}{2}-1}\choose{j-1}}$. The number of free crossings is $f=2m-2j-2k+2\ge 0$, which implies that $j\le m+1-k$ hence $j\le \min\{k,m+1-k\}$. There are ${{m-j-k+1+\frac{k+j-1}{2}-1}\choose{\frac{k+j-1}{2}-1}}={{m-\frac{k+j+1}{2}}\choose{\frac{k+j-1}{2}-1}}$ ways to distribute the free crossings symmetrically among the $k+j-1$ (symmetric) medium Seifert circles. The result follows since $\frac{k+j-1}{2}=\lfloor\frac{k+j}{2}\rfloor$ and $\frac{k+j+1}{2}=\lceil\frac{k+j}{2}\rceil$
when $k+j$ is odd. 

\medskip
Similarly, if $k-j$ is even, then we will treat the medium Seifert circle in the middle of the template as two. This creates $2j$ slots and the $k-j$ insertions are to be symmetrically distributed among them. See Figure \ref{Figure9} for a case of $s=2k=8$ and $2j=4$.
\begin{figure}[!hbt]
\begin{center}
\includegraphics[scale=0.6]{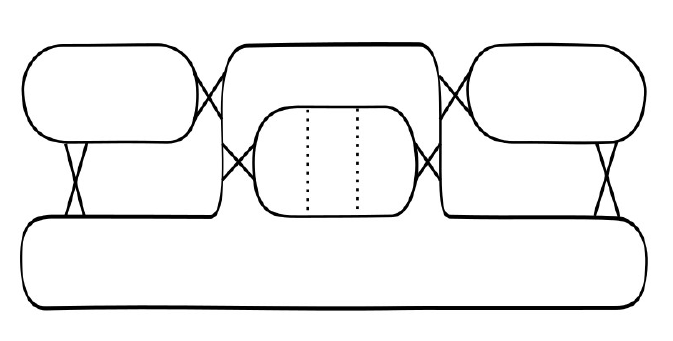} \qquad \includegraphics[scale=0.6]{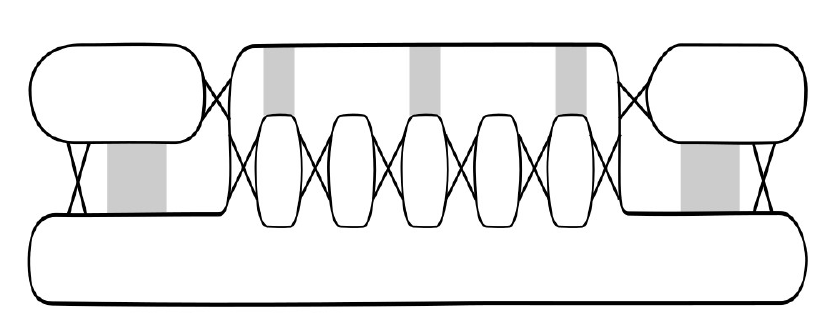}
\end{center}
\caption{A symmetric insertion distribution for an R-decomposition with 8 Seifert circles (right) built on an R-template with 4 Seifert circles (left).
\label{Figure9}}
\end{figure}

\medskip
The total number of such insertions is ${{\frac{k-j}{2}+j-1}\choose{j-1}}={{\frac{k+j}{2}-1}\choose{j-1}}$. In any of the resulting symmetric R-decompositions, there are $2j+k-j=k+j$ slots for the $2m-2j-2k+2$ free crossings to be  symmetrically distributed.
There are ${{m-j-k+1+\frac{k+j}{2}-1}\choose{\frac{k+j}{2}-1}}={{m-\frac{k+j}{2}}\choose{\frac{k+j}{2}-1}}$ ways to do so. The result follows since $\frac{k+j}{2}=\lfloor\frac{k+j}{2}\rfloor=\lceil\frac{k+j}{2}\rceil$
when $k+j$ is even. 

\medskip
Next let us consider the case when $n\ge 3$ is odd. In this case one free crossing must be placed right in the middle of the R-decomposition, which means that the Seifert circle in the middle must be a medium Seifert circle, so $j-k$ must be even. If $j-k$ is even, the counting formula is obtained similar to the case $n$ is even: we simply replace $n$ by $n-1$ since one crossing must be used in the middle. The result follows since $\frac{n-1}{2}=\lfloor \frac{n}{2}\rfloor$ when $n$ is odd.
\end{proof}

\section{The average genera of rational knots and links with a given crossing number}

Before such a formula can be derived, we need to know first which R-decompositions correspond to knots and which ones correspond to links. The following theorem is applicable to any oriented link diagram.

\begin{theorem}\label{T6}
Let $D$ be the diagram of an oriented link $L$ (it does not have to be a minimum diagram). Let $Cr(D)$  be the number of crossings in $D$, $s(D)$ be the number of Seifert circles in the Seifert circle decomposition of $D$ and $\mu(L)$ be the number of components of $L$, then $Cr(D)-s(D)=\mu(L)\ {\rm mod}2$.
\end{theorem}

\medskip
\begin{proof}
This can be proved by induction on the number of crossings in $D$. If $D$ contains no crossings, then it contains $\mu(L)$ topological circles that do not intersect each other, each such circle is a Seifert circle. So the statement of the theorem is trivially true. Assume that the statement is true for any diagram $D$ with $Cr(D)\le n$ where $n\ge 0$. Then consider the case when $D$ contains $n+1$ crossings. Let $D^\p$ be the diagram obtained from $D$ by smoothing a crossing of $D$. $\mu(D^\p)=\mu(D)+1$ if the strands at the crossing being smoothed belong to the same component, otherwise $\mu(D^\p)=\mu(D)-1$. Thus $\mu(D)-\mu(D^\p)=1\ {\rm mod}2$. On the other hand, $s(D)=s(D^\p)$, $Cr(D)=Cr(D^\p)+1$, hence $Cr(D)-s(D)=Cr(D^\p)-s(D^\p)+1=\mu(D^\p)+1=\mu(D)\ {\rm mod}2$ by the induction hypothesis.
\end{proof}

The following corollary is an immediate consequence of Theorem \ref{T6} since in the case that $D$ is the diagram of an oriented rational link, it can only have one or two components.

\begin{corollary}\label{C6}
An R-decomposition with $n$ crossings and $s$ Seifert circles is a knot if and only if $n-s\equiv 1\ {\rm mod}2$, and it is a link if and only if $n-s\equiv 0\ {\rm mod}2$.
\end{corollary}

\medskip
By \cite[Theorem 1]{Ernst1987}, the total number of rational knots with $n$ crossings (denoted by $RK_n$) is given by 
\begin{equation}\label{E_RKn}
RK_n=\frac{1}{3}\left(2^{n-2}+\frac{(-1)^{n+1}+1}{2}2^{\frac{n-1}{2}}+\frac{(-1)^{n+1}-1}{2}+1-(-1)^{\lfloor\frac{n}{2}\rfloor n}\right).
\end{equation}
On the other hand, by \cite[Theorem 4.3]{Diao2022}, the total number of oriented rational links with crossing number $n$ (denoted by $LK_n$) is given by
\begin{equation}\label{E_TLn}
RL_n=\frac{1}{3}\left(2^{n-2}+2(-1)^{n}\right)+\frac{(-1)^{n}+1}{2}2^{\frac{n-2}{2}}.
\end{equation}
Combining these formulas with Theorems \ref{T4} and \ref{T5}, we then obtain the main result of this paper in the following theorem.

\medskip
\begin{theorem} The average genus of rational knots with crossing number $n$, denoted by $\langle g\rangle_{K_n}$, is given by 
\begin{equation}\label{Kav}
\langle g\rangle_{K_n}=
\frac{1}{RK_n}\sum_{g=1}^{\lfloor\frac{n-1}{2}\rfloor} g\Psi_n(g),
\end{equation}
where $\Psi_n(g)=|R^I_n(n+1-2g)\cup R^{III}_n(n+1-2g)|+|RS^{III}_n(n+1-2g)|$ is the number of rational knots with crossing number $n$ and genus $g$, which can be expressed in the following explicit formula 
\begin{eqnarray}
&&
\sum_{j=1}^{\min\{g, \lceil\frac{n}{2}\rceil-g\}}{{\lfloor\frac{n}{2}\rfloor+j-g-1}\choose{2j+\frac{(-1)^n-1}{2}-1}}{ {\lceil\frac{n}{2}\rceil +g-j-1}\choose{\lfloor\frac{n}{2}\rfloor+j-g-1}}\label{E_RKng}\\
&+&
\left(\frac{1-(-1)^n}{2}\right)\sum_{j=1}^{\min\{g,\lceil\frac{n}{2}\rceil-g\}}\left(\frac{(-1)^{\lceil\frac{n}{2}\rceil+j-g}+1}{2}\right)
{{\lfloor\frac{j+\lceil\frac{n}{2}\rceil-g}{2}\rfloor-1}\choose{j-1}}{{\lfloor\frac{n}{2}\rfloor-\lceil \frac{j+\lceil\frac{n}{2}\rceil-g}{2}\rceil}\choose{\lfloor\frac{j+\lceil\frac{n}{2}\rceil-g}{2}\rfloor-1}}.\nonumber
\end{eqnarray}
Similarly, the average genus of oriented rational links with crossing number $n$, denoted by $\langle g\rangle_{L_n}$, is given by 
\begin{equation}\label{Lav}
\langle g\rangle_{L_n}=
\frac{1}{RL_n}\sum_{g=0}^{\lfloor\frac{n}{2}\rfloor} g\Phi_n(g),
\end{equation}
where $\Phi_n(g)=|R^I_n(n-2g)\cup R^{III}_n(n-2g)|+|RS^{III}_n(n-2g)|$ is the number of oriented rational links with crossing number $n$ and genus $g$, which has the following explicit formula 
\begin{eqnarray}
&&
\sum_{j=1}^{\min\{g+\frac{(-1)^n+1}{2}, \lfloor\frac{n}{2}\rfloor-g\}}{{j+\lfloor\frac{n-1}{2}\rfloor-g-1}\choose{2j-\frac{(-1)^n+1}{2}-1}}{ {\lfloor\frac{n}{2}\rfloor +g-j}\choose{j+\lfloor\frac{n-1}{2}\rfloor-g-1}}\label{E_RLng}\\
&+&
\left(\frac{1+(-1)^n}{2}\right)\sum_{j=1}^{\min\{g+1,\lfloor\frac{n}{2}\rfloor-g\}}
{{\lfloor\frac{\lfloor\frac{n}{2}\rfloor+j-g}{2}\rfloor-1}\choose{j-1}}{{\lfloor\frac{n}{2}\rfloor-\lceil \frac{\lfloor\frac{n}{2}\rfloor+j-g}{2}\rceil}\choose{\lfloor\frac{\lfloor\frac{n}{2}\rfloor+j-g}{2}\rfloor-1}}.\nonumber
\end{eqnarray}
\end{theorem}

\medskip
Notice that $g\ge 1$ in (\ref{Kav}) and $g\ge 0$ in (\ref{Lav}) since the minimum genus for any nontrivial knot is at least one, while it is possible that the minimum genus for a nontrivial link to be 0. The explicit formulas for $\Psi_n(g)$ and $\Phi_n(g)$ follow from Theorems \ref{T4} and \ref{T5}.
The plots in Figure \ref{K_plot} are the plots of $\langle g\rangle_{K_n}$ and $\langle g\rangle_{L_n}$ respectively, both are nearly linear.

\medskip
\begin{figure}[!ht]
\begin{center}
\begin{tikzpicture}
\begin{axis}[
    xmin=0, xmax=55,
    ymin=0, ymax=14,
    xtick={0,10,20,30,40,50},
    ytick={0,2,4,6,8,10,12},
    legend pos=north west,
    legend style={at={(0.5,-0.1)},anchor=north},
    ymajorgrids=true,
    grid style=dashed,
]
\addplot [
    domain=0:60, 
    samples=100, 
    color=red,
]
{0.249541338*x+0.100479927};
\addplot[
    color=black,
    mark=o,
    ]
     table{genus_knot.dat};
    \addlegendentry{$\langle g\rangle_{K_n}=0.2495n+0.1005$}
    \addlegendentry{$R^2=0.9999$}    
\end{axis}
\end{tikzpicture}
\hspace{1cm}
\begin{tikzpicture}
\begin{axis}[
    xmin=0, xmax=55,
    ymin=0, ymax=14,
    xtick={0,10,20,30,40,50},
    ytick={0,2,4,6,8,10,12},
    legend pos=north west,
    legend style={at={(0.5,-0.1)},anchor=north},
    ymajorgrids=true,
    grid style=dashed,
]

\addplot [
    domain=0:60, 
    samples=100, 
    color=red,
]
{0.249967706*x-0.41563816};

\addplot[
    color=black,
    mark=o,
    ]
     table{genus_link.dat};
    \addlegendentry{$\langle g\rangle_{L_n}=0.2499n-0.4156$} 
    \addlegendentry{$R^2=0.9999$}  
\end{axis}
\end{tikzpicture}

\end{center}
\caption{Left: The plot of $\langle g\rangle_{K_n}$, where the horizontal axis is the crossing number $n$ and the vertical axis is $\langle g\rangle_{K_n}$. Right: The plot of $\langle g\rangle_{L_n}$, where the horizontal axis is the crossing number $n$ and the vertical axis is $\langle g\rangle_{L_n}$. }\label{K_plot}
\end{figure}

\section{Further discussions}

Since $g=(n-s-\mu+2)/2< n/2$, we have $\langle g\rangle_{K_n}<n/2$ and $\langle g\rangle_{L_n}<n/2$. We propose the following conjecture based on the strong numerical  evidence as shown in Figure \ref{K_plot}.

\begin{conjecture}\label{bound}
$\lim_{n\to \infty}\frac{\langle g\rangle_{K_n}}{n}=\lim_{n\to \infty}\frac{\langle g\rangle_{L_n}}{n}=\frac{1}{4}$.
\end{conjecture}

\medskip
For each fixed number $s$ of Seifert circles, the numbers of rational knots or rational links form a sequence of positive integers with the crossing number as the running index. Some of these sequences appear in the Online Encyclopedia of Integer Sequences~\cite{OEIS}, and some can serve as generalizations of some existing sequences in the data base. 
Table \ref{t1}  lists the sequences $|R_n^I(s)\cup R_n^{III}(s)|$.  We note that the sequence $|R_{n+3}^I(4)\cup R_{n+3}^{III}(4)|$ is the sequence A000124 in the Online Encyclopedia of Integer Sequences~\cite{OEIS}, while the sequences $|R_{n+5}^I(5)\cup R_{n+5}^{III}(5)|$ and $|R_{n+7}^I(7)\cup R_{n+7}^{III}(7)|$ are the sequences A003600 and A257890 OEIS ~\cite{OEIS}. The sequences $|R_n^I(2k-1)\cup R_n^{III}(2k-1)|$ and $|R_n^I(2k)\cup R_n^{III}(2k)|$ ($k\ge 2$, $n\ge 2k$) are related by the equation 
$$
|R_n^I(2k-1)\cup R_n^{III}(2k-1)|+|R_n^I(2k)\cup R_n^{III}(2k)|=|R_{n+1}^I(2k)\cup R_{n+1}^{III}(2k)|.
$$
Also, the sequence $|R_{n+1}^I\cup R^{III}_{n+1}|$ is the Jacobsthal sequence (A001045) in the OEIS.

\begin{longtable}{|c|c|c|c|c|c|c|c|c|c|c|c|c|c|c|}
\caption{Each number in  the $(c(D),s(D))$ position represents the total number of rational knots/links with $s(D)$ Seifert circles and crossing number $c(D)$. The numbers in bold indicate links.}
\label{t1}\\
\hline
\backslashbox{$c(D)$}{$s(D)$}&2&3&4&5&6&7&8&9&10&11&12&13&14&$|R_n^I\cup R^{III}_n|$\\\hline\hline
2&\textbf{1}&&&&&&&&&&&&&1\\\hline
3&1&&&&&&&&&&&&&1\\\hline
4&\textbf{1}&1&\textbf{1}&&&&&&&&&&&3\\\hline
5&1&\textbf{2}&2&&&&&&&&&&&5\\\hline
6&\textbf{1}&3&\textbf{4}&2&\textbf{1}&&&&&&&&&11\\\hline
7&1&\textbf{4}&7&\textbf{6}&3&&&&&&&&&21\\\hline
8&\textbf{1}&5&\textbf{11}&13&\textbf{9}&3&\textbf{1}&&&&&&&43\\\hline
9&1&\textbf{6}&16&\textbf{24}&22&\textbf{12}&4&&&&&&&85\\\hline
10&\textbf{1}&7&\textbf{22}&40&\textbf{46}&34&\textbf{16}&4&\textbf{1}&&&&&171\\\hline
11&1&\textbf{8}&29&\textbf{62}&86&\textbf{80}&50&\textbf{20}&5&&&&&341\\\hline
12&\textbf{1}&9&\textbf{37}&91&\textbf{148}&166&\textbf{130}&70&\textbf{25}&5&\textbf{1}&&&683\\\hline
13&1&\textbf{10}&46&\textbf{128}&239&\textbf{314}&296&\textbf{200}&95&\textbf{30}&6&&&1365\\\hline
14&\textbf{1}&11&\textbf{56}&174&\textbf{367}&553&\textbf{610}&496&\textbf{295}&125&\textbf{36}&6&\textbf{1}&2731\\\hline
15&1&\textbf{12}&67&\textbf{230}&541&\textbf{920}&1163&\textbf{1106}&791&\textbf{420}&161&\textbf{42}&7&5461\\\hline
\end{longtable}

\medskip
Table \ref{t2}  lists the sequences $|RS_n^{III}(s)|$.  The sequence $|RS_n^{III}(4)|$ is the sequence A028242 in OEIS~\cite{OEIS}, while for the sequence $RS^{III}_n$ we observe that $|RS_{2n}^{III}|=2^n$ ($n\ge 1$) with $|RS_{2n}^{III}(2k)|=2^n$, $1\le k\le n$, being the coefficients of the binomial expansion of $(x+y)^n$. Furthermore, we observe that $|RS_{2n+1}^{III}|=|R_{n+1}^I\cup R^{III}_{n+1}|$ ($n\ge 1$) is also the Jacobsthal sequence (A001045) in the OEIS. 

\begin{longtable}{|c|c|c|c|c|c|c|c|c|c|c|c|}
\caption{Each number in  the $(c(D),s(D))$ position represents the total number of rational knots/links with $s(D)$ Seifert circles and crossing number $c(D)$ that have symmetric R-decompositions. The numbers in bold indicate links.}\label{t2}\\
  \hline
 \multicolumn{1}{|c|}{{\backslashbox{$c(D)$}{$s(D)$}}} &
  \multicolumn{1}{c|}{2} &
   \multicolumn{1}{c|}{4} &
   \multicolumn{1}{c|}{6} &
   \multicolumn{1}{c|}{8} &
   \multicolumn{1}{c|}{10} &
   \multicolumn{1}{c|}{12} &
   \multicolumn{1}{c|}{14} &
   \multicolumn{1}{c|}{16} &
   \multicolumn{1}{c|}{18} &
   \multicolumn{1}{c|}{20} &
   \multicolumn{1}{c|}{$|RS_n^{III}|$} \\
   \hline
\endfirsthead
  \hline
 \multicolumn{1}{|c|}{{\backslashbox{$c(D)$}{$s(D)$}}} &
  \multicolumn{1}{c|}{2} &
   \multicolumn{1}{c|}{4} &
   \multicolumn{1}{c|}{6} &
   \multicolumn{1}{c|}{8} &
   \multicolumn{1}{c|}{10} &
   \multicolumn{1}{c|}{12} &
   \multicolumn{1}{c|}{14} &
   \multicolumn{1}{c|}{16} &
   \multicolumn{1}{c|}{18} &
   \multicolumn{1}{c|}{20} &
   \multicolumn{1}{c|}{$|RS_n^{III}|$} \\
   \hline
\endhead
2&\textbf{1}&&&&&&&&&&1\\\hline
3&1&&&&&&&&&&1\\\hline
4&\textbf{1}&\textbf{1}&&&&&&&&&2\\\hline
5&1&{0}&&&&&&&&&1\\\hline
6&\textbf{1}&\textbf{2}&\textbf{1}&&&&&&&&4\\\hline
7&1&{1}&1&&&&&&&&3\\\hline
8&\textbf{1}&\textbf{3}&\textbf{3}&\textbf{1}&&&&&&&8\\\hline
9&1&{2}&2&{0}&&&&&&&5\\\hline
10&\textbf{1}&\textbf{4}&\textbf{6}&\textbf{4}&\textbf{1}&&&&&&16\\\hline
11&1&3&4&{2}&1&&&&&&11\\\hline
12&\textbf{1}&\textbf{5}&\textbf{10}&\textbf{10}&\textbf{5}&\textbf{1}&&&&&32\\\hline
13&1&{4}&7&{6}&3&{0}&&&&&21\\\hline
14&\textbf{1}&\textbf{6}&\textbf{15}&\textbf{20}&\textbf{15}&\textbf{6}&\textbf{1}&&&&64\\\hline
15&1&{5}&11&{13}&9&{3}&1&&&&43\\\hline
16&\textbf{1}&\textbf{7}&\textbf{21}&\textbf{35}&\textbf{35}&\textbf{21}&\textbf{7}&\textbf{1}&&&128\\\hline
17&1&{6}&16&{24}&22&{12}&4&{0}&&&85\\\hline
18&\textbf{1}&\textbf{8}&\textbf{28}&\textbf{56}&\textbf{70}&\textbf{56}&\textbf{28}&\textbf{8}&\textbf{1}&&256\\\hline
19&1&{7}&22&{40}&46&{34}&16&{4}&1&&171\\\hline
20&\textbf{1}&\textbf{9}&\textbf{36}&\textbf{84}&\textbf{126}&\textbf{126}&\textbf{84}&\textbf{36}&\textbf{9}&\textbf{1}&512\\\hline
21&1&{8}&29&{62}&86&{80}&50&{20}&5&{0}&341\\
\hline
\end{longtable}

Finally, Table \ref{t3} lists the sequences $R_n(s)$. We observe that for each fixed $n$, the distribution of the numbers $R_n(2)$, $R_n(3)$, ..., $R_n(n)$ has a bell shape with the maximum appearing at $s=(n+1)/2$ when $n$ is odd, and at $s=(n+2)/2$ when $n$ is even. Assuming that this is generally true, then it is not hard to see that $\langle g\rangle_{K_n}$ and $\langle g\rangle_{L_n}$ must grow at least linearly in terms of $n$.

\begin{longtable}{|c|c|c|c|c|c|c|c|c|c|c|c|c|c|c|}
\caption{The total numbers of rational knots and links with a given crossing number whose minimum projections have the same number of Seifert circles. Again, the numbers in bold indicate links.}
\label{t3}\\\hline
\backslashbox{$c(D)$}{$s(D)$}&2&3&4&5&6&7&8&9&10&11&12&13&14&$|\Lambda_{c(D)}|$\\\hline\hline
2&\textbf{2}&&&&&&&&&&&&&2\\\hline
3&2&&&&&&&&&&&&&2\\\hline
4&\textbf{2}&1&\textbf{2}&&&&&&&&&&&5\\\hline
5&2&\textbf{2}&2&&&&&&&&&&&6\\\hline
6&\textbf{2}&3&\textbf{6}&2&\textbf{2}&&&&&&&&&15\\\hline
7&2&\textbf{4}&8&\textbf{6}&4&&&&&&&&&24\\\hline
8&\textbf{2}&5&\textbf{14}&13&\textbf{12}&3&\textbf{2}&&&&&&&51\\\hline
9&2&\textbf{6}&18&\textbf{24}&24&\textbf{12}&4&&&&&&&90\\\hline
10&\textbf{2}&7&\textbf{26}&40&\textbf{52}&34&\textbf{20}&4&\textbf{2}&&&&&187\\\hline
11&2&\textbf{8}&32&\textbf{62}&90&\textbf{80}&52&\textbf{20}&6&&&&&352\\\hline
12&\textbf{2}&9&\textbf{42}&91&\textbf{158}&166&\textbf{140}&70&\textbf{30}&5&\textbf{2}&&&715\\\hline
13&2&\textbf{10}&50&\textbf{128}&246&\textbf{314}&302&\textbf{200}&98&\textbf{30}&6&&&1386\\\hline
14&\textbf{2}&11&\textbf{62}&174&\textbf{382}&553&\textbf{630}&496&\textbf{310}&125&\textbf{42}&6&\textbf{2}&2795\\\hline
15&2&\textbf{12}&72&\textbf{230}&552&\textbf{920}&1176&\textbf{1106}&800&\textbf{420}&164&\textbf{42}&8&5504\\\hline
\end{longtable}


\medskip
\begin{thebibliography}{99}
\bibitem{A} C.~Adams, {\em The Knot Book}, American Mathematical Soc., 1994.

\bibitem{BZ} G.~ Burde, H.~ Zieschang and M.~ Heusener
{\em Knots}, De Gruyter Studies in Mathematics \textbf{5},  2013.

\bibitem{Cohen} M.~Cohen, A Lower Bound on the Average Genus of a 2-Bridge Knot, preprint. \url{https://arxiv.org/pdf/2108.00563.pdf}

\bibitem{Diao2022} Y.~Diao, M. ~Finney, and D.~Ray, The Number of Rational Links with a Given Deficiency, {\em Journal of Knot Theory and its Ramifications} \textbf{30} 9 (2022), Article 2150065.

\bibitem{Gabai1986} D. ~Gabai, Genera of the Alternating Links, {\em Duke Mathematical Journal} \textbf{53} 3 (1986), 677--681.

\bibitem{DiaoErnst2005} Y.~Diao and C.~Ernst, The Growth Rate of Some Deficiency Zero Knot Classes, {\em International Journal of Pure and Applied Mathematics} \textbf{23} 4 (2005), 437--450.

\bibitem{DL}
  Y.~ Diao, C.~ Ernst, G.~ Hetyei and P.~ Liu,
A diagrammatic approach for determining the braid index of alternating
links, {\em Journal of Knot Theory and its Ramifications} \textbf{30} (5), 2150035 (2021). 

\bibitem{Ernst1987} 
C.~Ernst and D.~Sumners, The growth of the number of prime knots. {\em Mathematical Proceedings of the Cambridge Philosophical Society}, \textbf{102}(2) (1987), 303--315. 

\bibitem{OEIS}
OEIS Foundation Inc.\ (2011), {\em The On-Line Encyclopedia of Integer Sequences,}
published electronically at \url{http://oeis.org}. 
\end{thebibliography}
\end{document}